\documentclass[11pt]{article}
\usepackage{authblk}

\usepackage[utf8]{inputenc}
\usepackage[T1]{fontenc}
\usepackage{mathpazo}

\RequirePackage[english]{babel}

\RequirePackage[a4paper,margin=2.5cm]{geometry}
\RequirePackage{parskip}

\makeatletter
\def\@fnsymbol#1{\ensuremath{\ifcase#1\or *\or \dagger\or \ddagger\or \mathsection\or \|\or **\or \dagger\dagger \or \ddagger\ddagger \else\@ctrerr\fi}}
\makeatother

\usepackage{xcolor}
\definecolor{cblue}{RGB}{0,70,140}
\definecolor{cgreen}{RGB}{100,140,0}
\definecolor{cred}{RGB}{190,10,50}

\usepackage[shortlabels]{enumitem}
\setlist{itemsep=0ex,topsep=0ex,parsep=0.4ex}

\usepackage{amsmath,amsfonts,amssymb,amsthm,mathtools,thmtools,thm-restate,bbm,centernot}

\usepackage[hyphens]{url} 
\usepackage[hidelinks,colorlinks,pagebackref]{hyperref} 
\hypersetup{citecolor=cgreen,linkcolor=cblue,urlcolor=cblue}
\usepackage[capitalise,nameinlink,noabbrev,compress]{cleveref} 

\renewcommand*{\backref}[1]{}
\renewcommand*{\backrefalt}[4]{
	\ifcase #1 Not cited.%
	\or $\uparrow$#2%
	\else $\uparrow$#2%
	\fi%
}

\let\oldbibliography\bibliography
\renewcommand{\bibliography}[1]{
  {

    \hypersetup{linkcolor=cred}
    \bibliographystyle{bibstyle}
    \oldbibliography{#1}
  }
}

\theoremstyle{plain}
\newtheorem{theorem}{Theorem}[section]
\newtheorem{lemma}[theorem]{Lemma}
\newtheorem{claim}[theorem]{Claim}

\newtheorem{problem}[theorem]{Problem}

\theoremstyle{definition}

\makeatletter
\renewenvironment{proof}[1][\proofname]
{\par\pushQED{\qed}
	\normalfont\topsep6\p@\@plus6\p@\relax\trivlist
	\item[\hskip\labelsep\bfseries#1\@addpunct{.}]
	\ignorespaces}
{\popQED\endtrivlist\@endpefalse}
\makeatother

\def\cro{{\mbox {\sc cr}}}
\def\ocr{{\mbox {\sc ocr}}}
\def\pcr{{\mbox {\sc pcr}}}

\title{Improvement of the Crossing Lemma for the odd-crossing number}

\author{N\'ora Alm\'asi\thanks{Supported by National Research, Development and Innovation Office, NKFIH ADVANCED-152590. Email:~\href{mailto:almasi.nora@cs.bme.hu}{\tt almasi.nora@cs.bme.hu}.} }
\author{G\'eza T\'oth\thanks{Supported by National Research, Development and Innovation Office, NKFIH
ADVANCED-152590 and ERC Advanced Grant “GeoScape" 882971. Email:~\href{mailto:geza@renyi.hu}{\tt geza@renyi.hu}.}}

\affil{HUN-REN Alfr\'ed R\'enyi Institute of Mathematics,
Budapest, Hungary}
\affil{Department of Computer Science and Information Theory \\
Budapest University of Technology and Economics, Budapest, Hungary}

\begin{document}

\maketitle

 \begin{abstract}
A graph is called $1$-planar if it can be drawn in the plane 
so that every edge is crossed by at most one other edge.
It is known that the maximum number of edges of $1$-planar graphs
is $4n-8$. We prove the following generalization.
A graph is $1$-odd-planar if it can be drawn in the plane 
so that every edge is crossed {\em an odd number of times}
by at most one other edge.
We show that the maximum number of edges of $1$-odd-planar graphs
is also $4n-8$.
As a consequence, we improve the constant coefficient in the Crossing Lemma for the odd-crossing number. 
\end{abstract}

\section{Introduction}

Throughout the paper, all graphs are assumed to be simple, that is, they contain neither loops nor parallel edges.
In a  {\em drawing} ${\cal D}$
of a graph $G$, vertices are represented by (or mapped to) distinct points in the plane, 
edges are represented by simple curves connecting the corresponding points.
If it is clear from the context, these points and curves are also called vertices and edges, respectively. 
We assume that an edge does not go through a vertex (apart from its endpoints) and no three edges cross at a common point. 
The {\em rotation system} of $G$ in drawing ${\cal D}$ 
is the set of cyclic orders of incident edges at each vertex.

The {\em crossing number} $\cro(G)$ of a graph $G$ is the minimum number of edge crossings (or crossing points)
over all drawings of $G$. 
The Crossing Lemma, proved  by Ajtai, Chv\'atal, Newborn, Szemer\'edi
\cite{ACNS82} and by Leighton \cite{L83}, is undoubtedly the most important inequality for crossing numbers.

\medskip

\noindent {\bf Crossing Lemma} {\em If a simple graph $G$ of $n$ vertices has $m\ge 4.5n$
edges, then $\cro(G)\ge \frac{1}{60.75}\frac{m^3}{n^2}$ edges.}

\medskip

The bound is tight, apart from the value of the constant 
\cite{PRTT06}, \cite{CSSW20}.

The constant above follows from the beautiful probabilistic argument of
Chazelle, Sharir and Welzl \cite{AZ99}. 
The constant has been  improved
in  five steps \cite{PT97}, \cite{PRTT06}, \cite{A19}, 
\cite{BK26}, 
the best bound is due to 
Büngener et al. 
\cite{BFKP26}, namely 
$\cro(G)\ge \frac{1}{27.3}\frac{m^3}{n^2}$, when $m\ge 7.39n$.

There are many variants of the crossing number \cite{S12}, \cite{S18}.
The most interesting and studied are the {\em pair-crossing number} 
$\pcr(G)$ and the {\em odd-crossing number}
$\ocr(G)$.
The {\em pair-crossing number} $\pcr(G)$ is the minimum number 
of pairs of crossing edges
over all drawings of $G$.
The {\em odd-crossing number} $\ocr(G)$ is the minimum number 
of pairs of edges that cross an odd number of times,
over all drawings of $G$.
Clearly, every graph $G$ satisfies
$\ocr(G)\le\pcr(G)\le\cro(G)$.

In \cite{PT00b}, some further variants were introduced in order to study the role of crossings between adjacent edges.
For each of $\cro$, $\pcr$, $\ocr$, they introduced three counting rules:

\smallskip

\noindent {\bf Rule $+$:} Only those drawings are considered, where adjacent edges cannot cross.

\noindent {\bf Rule $0$:} Adjacent edges can cross and their crossings are counted as well.

\noindent {\bf Rule $-$:} Adjacent edges can cross and their crossings are not counted.

\smallskip

Combining these rules with the three crossing numbers yields nine possibilities.
However, it is easy to see that
$\cro_{+}=\cro$ \cite{PT00b}.
On the other hand, as observed in \cite{KT23}, 
regarding Rule $+$ for the odd-crossing number, it seems more natural to assume that adjacent edges cross an {\em even number of times} than
to assume that they do not cross at all.
So, let $\ocr_{*}(G)$ be the minimum number of odd-crossing pairs of edges over all drawings of $G$
where adjacent edges cross an even number of times. 
Therefore, we have nine versions, see Table 1, taken from \cite{KT23}. 
For a survey about crossings between adjacent edges, see \cite{F26}.
In this table, it follows directly from their definitions 
that values do not decrease if we move to the right or up.

\smallskip

{\scriptsize

\setlength{\unitlength}{1cm}
\begin{picture}(10, 3.5)
\put(3.5, 0.2){

\begin{picture}(9.5, 3)

\put(0, 0){\line(1, 0){8.5}}
\put(0, 1){\line(1, 0){8.5}}
\put(0, 2){\line(1, 0){6.5}}
\put(0, 3){\line(1, 0){8.5}}


\put(0, 0){\line(0, 1){3}}
\put(1.5, 0){\line(0, 1){3}}
\put(4.5, 0){\line(0, 1){3}}
\put(6.5, 0){\line(0, 1){3}}
\put(8.5, 0){\line(0, 1){3}}

\put(0.2, 0.3){Rule {--}}
\put(0.2, 1.3){Rule 0}
\put(0.2, 2.3){Rule +}


\put(7.1, 1.9){$\cro(G)$}
\put(7.1, 0.4){$\cro_{-}(G)$}

\put(4.6, 2.4){$\pcr_{+}(G)$}
\put(4.6, 1.4){$\pcr(G)$}
\put(4.6, 0.4){$\pcr_{-}(G)$}

\put(1.6, 2.4){$\ocr_{*}(G)\le \ocr_{+}(G)$}
\put(1.6, 1.4){$\ocr(G)$}
\put(1.6, 0.4){$\ocr_{-}(G)$}

\end{picture}}
\end{picture}

\medskip
}
\centerline{{\bf Table 1.} Nine versions of the crossing number.}

\medskip

Interestingly, the proof of the Crossing Lemma by 
Chazelle, Sharir and Welzl \cite{AZ99} works for all of these versions \cite{PT00a}, \cite{T26}
of the crossing numbers, but the improvements do not.

There are  only two improvements for other versions, by 
Ackerman and Schaefer \cite{AS14} and by Karl and T\'oth \cite{KT23}, respecively. 

\begin{theorem}[\cite{AS14}]\label{AckermanSchaefer}
Suppose that $G$ has $n$ vertices and $m\ge 6.75n$ edges, then 
$\pcr_{+}(G)\ge \frac{1}{34.2}\frac{m^3}{n^2}$.
\end{theorem}

\begin{theorem}[\cite{KT23}]\label{KarlToth}
Suppose that $G$ has $n$ vertices and $m\ge 6n$ edges, then 
$\ocr_{+}(G)\ge\ocr_{*}(G)$
$\ge \frac{1}{54}\frac{m^3}{n^2}$.
\end{theorem}

For the other versions, the best currently known constant is $\frac{1}{60.75}$.

In this note, we improve Theorem \ref{KarlToth}.

\begin{theorem}\label{oddcrlemma}
  Suppose that $G$ has $n$ vertices and $m\ge 5.25n$ edges. Then
  $\ocr_{+}(G)\ge\ocr_{*}(G)$
  $\ge \frac{1}{41.4}\frac{m^3}{n^2}$. 
\end{theorem}

All these improvements are based on the same idea.
The first step is to find many crossing pairs in sparse graphs.
Then this bound is applied to a random subgraph to get the general bound.
A graph $G$ is called $k$-planar if it can be drawn in the plane such that
there are at most $k$ crossings on each edge. Such a drawing is called a $k$-plane drawing.
Let $m_k(n)$ denote the maximum number of edges of a $k$-planar graph of $n$
vertices.
Clearly, $m_0(n)=3n-6$.
It is known that  $m_1(n)=4n-8$ for $n\ge 12$,
$m_2(n)\le 5n-10$ and it is  tight for infinitely
many values of $n$ \cite{PT97},
$m_3(n)\le 5.5n-11$, $m_4(n)\le 6n-12$,
which are tight up to an additive constant
\cite{PRTT06}, \cite{A19}.
These bounds, and some other observations are used for the above mentioned improvements for the constant in the Crossing Lemma for the classical crossing number $\cro(G)$.

Similarly, a graph $G$ is called $k$-odd-planar if it can be drawn in the plane such that
any edge is crossed {\em an odd number of times} by at most $k$ other edges
(edges crossing an even number of times do not count). Such a drawing is called a
 $k$-odd-plane drawing.

Let $m^{\mbox{\scriptsize odd}}_k(n)$ denote the maximum number of edges of a
$k$-odd-planar
graph with $n$ vertices.
Clearly, we have $m^{\mbox{\scriptsize odd}}_k(n)\ge m_k(n)$
and, by the weak Hanani-Tutte theorem
 \cite{H34}, \cite{PSS07}, we have $m^{\mbox{\scriptsize odd}}_0(n)=3n-6$.

Karl and T\'oth \cite{KT23} proved that $m^{\mbox{\scriptsize odd}}_1(n)\le 5n-9$, and used this result to show Theorem \ref{KarlToth}.
Here we improve this bound to $4n-8$, matching the classical upper bound for 1-planar graphs, that is for $m+1(n)$. In particular, allowing arbitrarily many even crossings between edges does not increase the maximum number of edges.

\begin{theorem}\label{thm:4n-8} 
Every simple 1-odd-planar graph with $n\geq 3$ vertices has at most $4n-8$ edges, i.e., $m_1^{\mathrm{odd}}(n)\leq 4n-8$.
Moreover, this bound is attained for every $n\geq 12$.

\end{theorem}

In Section \ref{sec:4n-8} we prove Theorem \ref{thm:4n-8}, and then, in Section \ref{sec:crlemma} we prove Theorem \ref{oddcrlemma}.

\section{On 1-odd-planar graphs}\label{sec:4n-8}

In a drawing of $G$, a subgraph $H\subseteq G$ is called an {\em even-crossing subgraph} if 
any two edges of $H$ cross an even number of times. 
If two edges cross an even (odd) number of times, then they are called an {\em even 
(odd) pair}.  
If $(e,f)$ is an odd pair in an 1-odd-plane drawing, we say that $e$ is the odd pair of $f$ (and $f$ is the odd pair of $e$). We use this terminology only in 1-odd-planar drawings where it really makes sense.
For any two edges $e$ and $f$, let $\cro(e, f)$ denote the number of crossings between $e$ and $f$.

Pelsmajer, Schaefer, and \v Stefankovi\'c \cite{PSS07} gave a very nice and simple proof of the weak Hanani-Tutte theorem
with a slightly 
stronger statement.

\begin{lemma}[Weak Hanani-Tutte Theorem \cite{PSS07}]\label{weakHTPSS}
If $G$ can be drawn in the plane so that any two edges cross an even number of times, then $G$ is planar. Moreover, it has a plane drawing with the same rotation system.

\end{lemma}

We prove the following generalization.

\begin{lemma}\label{redraw}
Let $G$ be a graph drawn in the plane and let $H$ be a connected spanning, 
even-crossing subgraph. 
We can redraw $G$ so that \\
(i) the rotation system is the same \\
(ii) odd and even  pairs of edges are the same as in the original drawing \\
(iii) edges of $H$ do not cross each other.
\end{lemma}

\begin{proof}
By the assumption, there is an even-crossing spanning tree $T$. 
We will contract all edges of $T$ recursively as follows.
Let $uv$ be an edge of $T$. Suppose that edges 
$e_1, \ldots, e_l$ 
cross $uv$  an odd number of times, 
and $f_1, \ldots f_m$ cross it an even number of times.
For each of $f_1, \ldots, f_m$ pull all crossings along the edge $uv$ through $u$, 
to all adjacent edges, without introducing new crossings among the $f_i$'s. 
Now they do not cross $uv$, and the parity of the number of crossings between any two edges remains unchanged.
Next, for each of  $e_1, \ldots, e_l$,  pull all crossings along the edge $uv$ through $u$, 
to all adjacent edges, {\em except one}, which is pulled exactly into $u$. 
Thus, for each $e_i$, exactly one crossing is pulled into $u$.
Now the edge $uv$ is crossing-free and the parity of the number of crossings between 
any two edges remained the same. But the curves representing $e_1, \ldots, e_l$ pass through $u$. 
Divide each of them into two edges at $u$. Contract the edge $uv$ to a vertex $w$.
Repeat this step for each edge of $T$. See Fig \ref{fig:osszehuzas}.

\begin{figure}[ht!]
    \centering
    \includegraphics[width=0.7\textwidth]{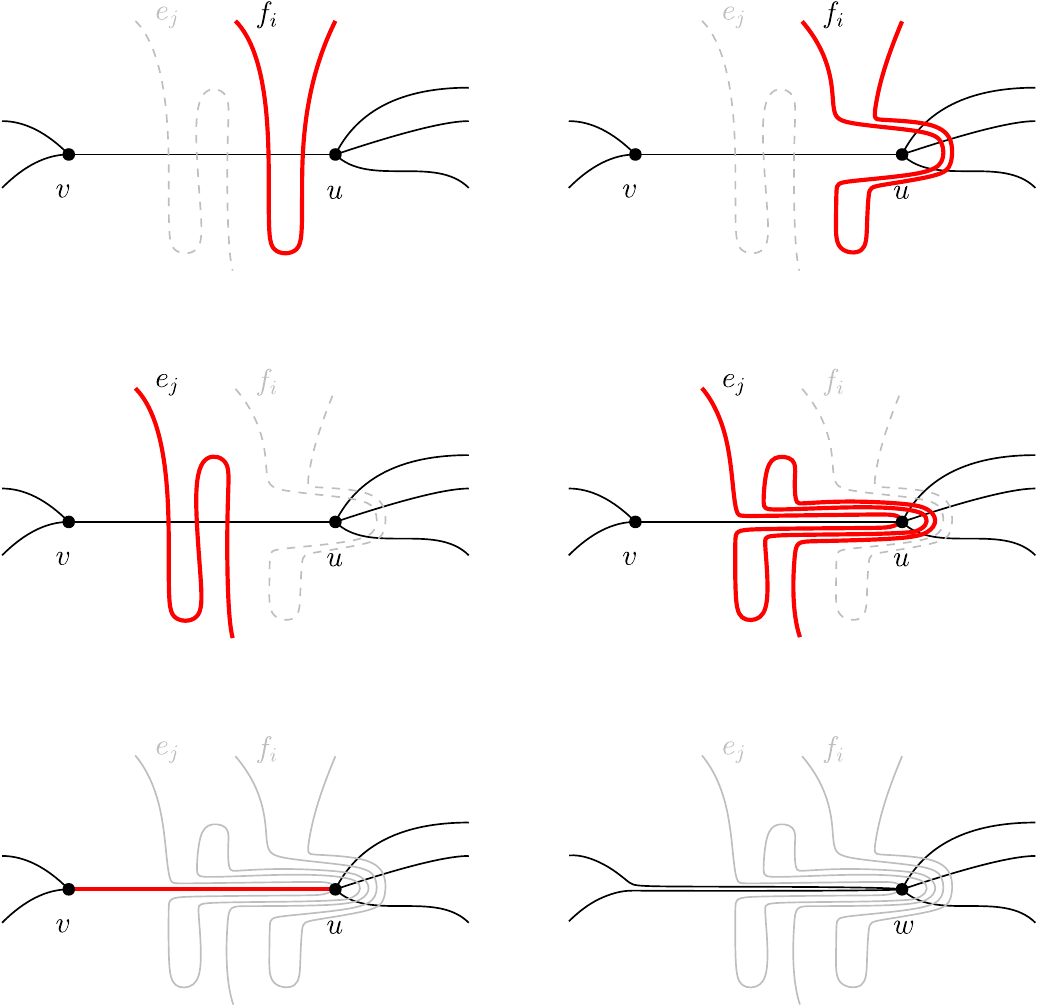}
    \caption{First pull all $f_i$ through $u$, then pull all $e_j$ as well, and finally contract the edge $uv$.}
    \label{fig:osszehuzas}
\end{figure}

Eventually we obtain a graph $G'$ with drawing ${\cal D}_2$,
which has only one vertex $v$ 
and all of its edges are loops.
Each loop corresponds to an edge, or a segment of an edge of $G$. 
Let $e$ and $f$ be two edges of $G\setminus T$ and suppose that $e$ corresponds to loops $e_1, \ldots, e_l$, 
$f$ corresponds to 
$f_1, \ldots, f_m$. Then, $\sum_{i=1}^l\sum_{j=1}^m\cro(e_i, f_j)\equiv \cro(e,f)$ $\bmod\  2$, 
that is, the total number of crossings between the pieces is the same mod $2$ 
as the number of crossings between the two original edges in ${\cal D}_1$.
Consider two loops, $e$ and $f$, and the rotation system at $v$. 
We either have $\ldots e\ldots e\ldots f\ldots f$
or $\ldots e\ldots f\ldots e\ldots f$. In the first case, $e$ and $f$ cross an even number of times, and in the second 
case they cross an odd number of times, regardless of the exact drawing of $e$ and $f$. 
We execute the following two operations as many times as possible.

If any two loops, say, $e$ and $f$, cross in two points,  
$p$ and $q$, then switch the segments of $e$ and $f$ between $p$ and $q$ and thus eliminate the two crossings.
If a loop intersects itself at $p$, then eliminate the 
self-crossing by a small modification at $p$. See Fig \ref{fig:removecr}.

\begin{figure}[ht!]
    \centering
    \includegraphics[width=0.6\textwidth]{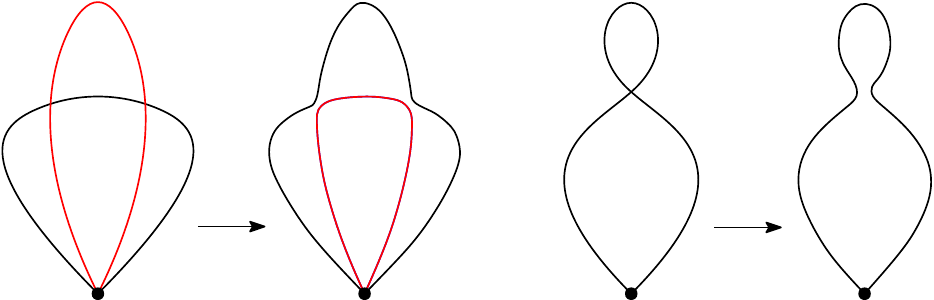}
    \caption{Loops cross at most once.}
    \label{fig:removecr}
\end{figure}


Both operations decrease the total number of crossings but do not change the rotation system at $v$. 
Therefore, after finitely many steps, we arrive at a situation where none of them 
can be executed, that is, any two loops cross at most once and there is no self-crossing. 
Moreover, the procedure does not change the parity of crossings between two edges since 
it does not change the rotation system. Let ${\cal D}_3$ be the resulting drawing of $G'$.

Observe that during the procedure in the proof of Lemma \ref{redraw}, those edges that cross the edges of $T$ an even number of times are not subdivided, so each of them corresponds to one loop of $G'$.

Now reverse the edge-contraction procedure, take contracted edges in the reverse order, and for each of them, do the following.  
Let $uv$ be the next edge, contracted to vertex $w$.
Replace $w$ by the two original vertices $u$ and $v$, 
very close to each other, connected by a crossing-free edge $uv$. 
Connect the edges at $w$ to their original vertex, $u$ or $v$. 
Suppose that in the contraction of $uv$, edge $e$ was subdivided into $e_1$ and $e_2$, both incident to $w$.
Then move the corresponding endpoints to a point $p$ on the edge $uv$ and merge  
$e_1$ and $e_2$ to get the original edge $e$. Then $p$ becomes a crossing between $uv$ and $e$.
We can assume that apart from these crossings of the last type, we do not create any additional crossings.

After completing these redrawing steps for all edges of $T$ in the reverse order, we obtain drawing ${\cal D}_4$ 
of $G$. In ${\cal D}_4$ the rotation system is the same as in ${\cal D}_1$. 
If two edges of $G$ in 
${\cal D}_1$ crossed an odd (even) number of times, then they also cross an odd (even) number of times in ${\cal D}_4$. 
Edges of $T$ do not cross any edge. 
Edges of $H\setminus T$ were never subdivided in the procedure,
so in ${\cal D}_3$  they do not cross. Therefore, they do not cross in ${\cal D}_4$ either.
This concludes the proof of the lemma. \end{proof}

\smallskip

\begin{proof}[Proof of Theorem \ref{thm:4n-8}]
Since $m_1(n)\le m_1^{\mathrm{odd}}(n)$ and $m_1(n)=4n-8$ for $n\ge 12$,   
we only have to prove that 
$m_1^{\mathrm{odd}}(n)\le 4n-8$. 
Suppose that $G$ is a $1$-odd-planar graph with $n$ vertices and $m$ edges. 
We may assume that $G$ is connected. Indeed, suppose that $G$ has $a$ isolated vertices, $b$ components with two vertices, and $c$ components with at least three vertices. If $c=0$, then $m\leq b\leq 4n-8$ holds trivially, since $n\geq 3$. Otherwise, applying the bound to each component with at least three vertices separately, we obtain
$m\leq b+4(n-a-2b)-8c\leq 4n-8$.
Thus, it suffices to prove the bound for connected graphs.


{\bf First, assume that $G$ has a $1$-odd-planar drawing with a 
connected spanning even-crossing subgraph.}

\smallskip

Among all $1$-odd-planar drawings of $G$,
let ${\cal D}_1$ be one where the maximal connected spanning even-crossing subgraph $H$ 
has the maximum number of edges. 
By the assumptions, every edge has at most one odd pair.

\smallskip

Let $\alpha$ be a (possibly self-crossing) closed curve. It divides the plane into connected regions. We can color them with black and white
such that two regions with a common boundary have different colors.
The union of the black (resp white) regions is called the two {\em sides} of the curve $\alpha$.

\smallskip

Return to the proof of Theorem \ref{thm:4n-8}.
Apply the procedure of Lemma \ref{redraw} and 
consider the drawing ${\cal D}_4$ of $G$. 
For simplicity, call the edges of $H$ {\em red}, the other edges {\em blue}.
By the assumption and Lemma \ref{redraw},
it has the following properties.
(i) the red edges do not cross each other
(ii) every blue edge crosses exactly one red edge an odd number of times
(iii) blue edges cross each other an even number of times.

Since the red edges do not cross each other, we have a plane drawing of $H$. 
Let $N(H)$, $E(H)$, $F(H)$ denote the numbers of vertices, edges, faces of $H$. 
Clearly, $N(H)=n$. By {\em face} we always mean a face of $H$.
A blue edge $uv$ {\em belongs to} a face $\Phi$ if its initial piece, very close to its endpoint, 
is in $\Phi$. One blue edge belongs to two faces, with multiplicity, that is, the two faces can be identical.
For any face $\Phi$, let $|\Phi|$ denote the number of edges on its boundary, counted with multiplicity. 
Let $\Phi_1, \ldots, \Phi_F$ denote  the faces where $F=F(H)$.
Let $E_b$ be the number of blue edges. We have to show that $E(G)=E(H)+E_b\le 4n-8$.

Let $\Delta$ be a triangular face. Let $u$, $v$, $w$ be its vertices in counterclockwise direction.

\begin{claim}\label{odd-belong}
A blue edge $e$ belongs to $\Delta$ if and only if it intersects one of its sides an odd number of times. \end{claim}

\begin{proof}
Suppose that a blue edge $e$ belongs to $\Delta$.  Let $u$ and $x$ be 
its endpoints.
Since the initial part of $e$, near $u$ is inside $\Delta$ and its other endpoint, $x$ is outside, it intersects one of its sides, say $f$, an odd number of times. 
Suppose now that $e$ intersects a side of $\Delta$, $f$ an odd number of times. Since it intersects the other two sides an even number of times, one of the initial parts is inside, 
the other one is outside $\Delta$. So $e$ belongs to $\Delta$.
\end{proof}

\begin{claim}\label{nemszomszedos}
Suppose that the blue edge $ux$ belongs to $\Delta$ and it intersects $uv$ an odd number of times.
Then in the cyclic order of edges at $u$, $ux$ and $uv$ are not consecutive. 
\end{claim}

\begin{proof}
Suppose they are consecutive. Switch their cyclic order in a very small neighborhood of $u$, by adding an additional crossing between them. Now $ux$ intersects all edges an even number of times. So, in this drawing, $H+ux$ is an even-crossing subgraph, contradicting the maximality of $H$. 
\end{proof}

\begin{claim}\label{otherface}
Suppose that edge $e$ belongs to $\Delta$ and its odd pair is side $f$ of 
$\Delta$. Let $\Phi$ be the face on the other side of $f$. Then $e$ also belongs to $\Phi$.
\end{claim}

\begin{proof}
Let $u$ and $x$ be the endpoints of $e$.
Since $e$ intersects $f$ an odd number of times, 
it intersects all other red edges an even number of times, so it intersects the boundary of $\Phi$ an odd number of times. Since the initial part of $e$ near $u$ is outside $\Phi$, the other initial part, near $x$, is inside. So, $e$ belongs to $\Phi$. 
\end{proof}

\begin{claim}\label{haromszog}

There is a face $\Phi$ with the property that
(i) all blue edges that belong to $\Delta$, also belong to $\Phi$, and (ii) 
if there are at least two such edges then $|\Phi|\ge 5$.

\end{claim}

\textbf{Proof.}
Clearly, no blue edge can belong to $\Delta$ twice since there are no loops or multiple edges in $G$.
We distinguish several cases.

{\bf 1.
There are exactly two blue edges that belong to $\Delta$.}

Up to symmetry, there are three different subcases distinguished 
by the incidence and cyclic order 
of these edges and their odd pairs, 
as shown in Figure \ref{fig:2-felel-2-csucs}.
All other scenarios are ruled out by Claim \ref{nemszomszedos}.

\begin{figure}[ht!]
    \centering
    \includegraphics[width=0.9\textwidth]{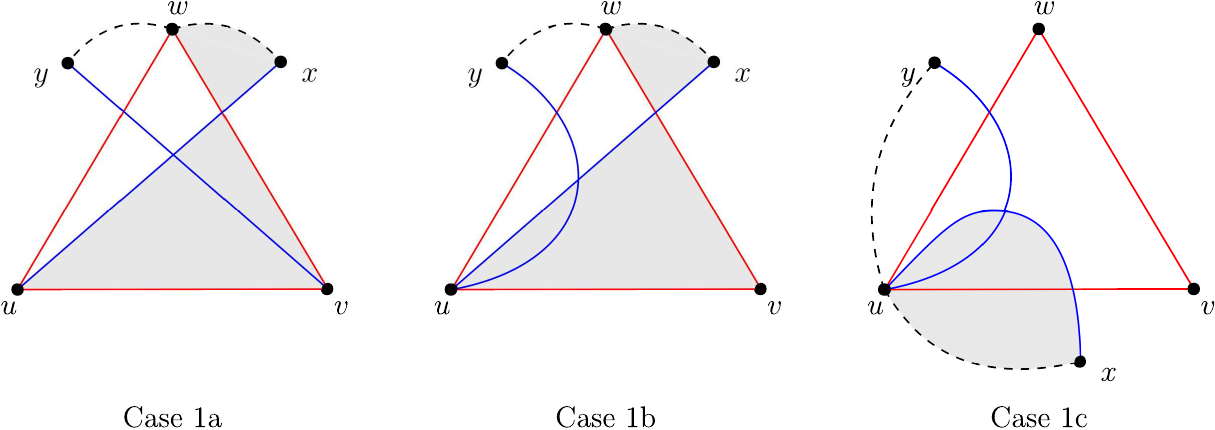}
    \caption{Two edges belong to $\Delta$.}
    \label{fig:2-felel-2-csucs}
\end{figure}



\textbf{1a. The two edges are incident to different vertices of $\Delta$.} 

Let $ux$ and $vy$ be the two edges.
By Claims \ref{nemszomszedos} and \ref{odd-belong}, the odd pair of $ux$ is $vw$ and the odd pair of $vy$ is $uw$.
Let $\Phi_x$ ($\Phi_y$) be the face on the other side of $vw$ ($uw$). We have to show that $\Phi_x=\Phi_y$. 
Suppose that they are different. 
Consider the closed, oriented curve $\gamma$ formed by the directed edges $\overrightarrow{xu}-\overrightarrow{uv}-\overrightarrow{vw}$ 
and the boundary of $\Phi_x$ from $w$ to $x$.
Let $\alpha$ be the boundary curve of $\Phi_y$ between $w$ and $y$, not containing $u$ and let $A$ and $B$ be the two sides of
$\gamma$.  At $v$, the initial segment of $vy$ lies on the {\em left} side of $\gamma$, suppose without loss of generality 
that it is side $A$. The edge $vw$ is crossed an even number of times by all edges, except $ux$, 
therefore, $vw$ contains an odd number of self-crossings of $\gamma$. Consequently, at $w$, $A$ is on the {\em right} side  
of $\gamma$ at $w$, while the initial part of $\alpha$ is on the right side, so it belongs to $B$.
Since every edge of $\gamma$ crosses every 
edge of $\alpha$ an even number of times, $y$ is also in side $B$. 
This is a contradiction, because $vy$ and $\gamma$ cross an even number of times.

\textbf{1b. The two edges are incident to the same vertex of $\Delta$, and the odd pair of one of them is the opposite side of $\Delta$.}

Let $ux$ and $uy$ be the two edges.
By Claims \ref{nemszomszedos} and \ref{odd-belong}, we may assume that the odd pair of $ux$ is $vw$, the odd pair of $uy$ is $uw$, and the cyclic order of edges at $u$ is $uw, ux, uy, uv$.
Let $\Phi_x$ ($\Phi_y$) be the face on the other side of $vw$ ($uw$). We have to show that $\Phi_x=\Phi_y$. Suppose that they are different. 
The argument is the same as the one in case \textbf{1a}.
We use the same notation and can again conclude that the initial parts of $uy$ lie on different sides of $\gamma$, while it crosses $\gamma$ an even number of times, which is again a contradiction.

\textbf{1c. The two edges are incident to the same vertex of $\Delta$, and the odd pairs of them are the two adjacent sides of $\Delta$.} 

Let $ux$ and $uy$ be the two edges.
By Claims \ref{nemszomszedos} and \ref{odd-belong}, we may assume that the odd pair of $ux$ is $uv$, the odd pair of $uy$ is $uw$, and the cyclic order of edges at $u$ is $uw, ux, uy, uv$.
Let $\Phi_x$ (resp. $\Phi_y$) be the face on the other side of $uv$ (resp. $uw$). 
We have to show that $\Phi_x=\Phi_y$. Suppose that they are different. 
Consider the closed curve $\gamma$ formed by the edge ${xu}$
and the boundary of $\Phi_x$ from $u$ to $x$.
Let $\alpha$ be the boundary curve of $\Phi_y$ between $u$ and $y$, 
not containing $w$. 
At $u$, the initial parts of $uy$ and $\gamma$ are on different sides of $\gamma$. 
Since $\alpha$ and $\gamma$ cross an even number of times, $y$ and the initial part of $uy$ at $u$ are also on 
different sides of $\gamma$. But it is a contradiction, since $uy$ and $\gamma$ also cross an even number of times.

{\bf 2. 
There are more than two blue edges that belong to $\Delta$.}

Clearly, at most three blue edges belong to $\Delta$.

\textbf{2a. The three edges are incident to at least two different vertices of $\Delta$.}

Let $ux$ be the only one out of the three edges that is incident to $u$ and let $e$ and $f$ be the other two. Suppose that 
$ux$ also belongs to $\Phi_x$. 
Then, by the argument in case \textbf{1a} for $ux$ and $e$, they 
belong to the same face, $\Phi_x$, and
similarly, $ux$ and $f$ also belong to the same face, $\Phi_x$.

\textbf{2b. All three edges are incident to the same vertex of $\Delta$.}

Suppose that all three edges are incident to $u$. 
Let the cyclic order of edges at $u$ be $uv$, $ux$, $uy$, $uz$, $uw$.
Each of the three edges 
has a different side of $\Delta$ as the odd pair.
If the odd pair of $uy$ is 
the opposite side, $vw$, then both pairs of $ux, uy$ and $uy, uz$ correspond to case \textbf{1b}, thus all three edges belong to the same face $\Phi$.
If the odd pair of $ux$ is the 
opposite side, $vw$, then the pair $ux, uz$ corresponds to case \textbf{1b}, and $uy, uz$ 
correspond to case \textbf{1c}, and thus all three edges belong to the same face $\Phi$.
If the odd pair of $uz$ is the 
opposite side, $vw$, then we are also done by symmetry. 

\smallskip

It remains to show that $| \Phi| \geq 5$.
Suppose not, and let $e$ and $f$ belong to $\Delta$ and $\Phi$. 
Then the odd pairs of $e$ and $f$ are common sides of $\Delta$ and $\Phi$.  
Therefore, $\Delta$ and $\Phi$ share at least two sides, say $uw$ and $wv$. Consequently, all vertices of $\Delta$, $u$, $v$, and $w$, are also vertices of $\Phi$. Consequently, $|\Phi|=4$, say, 
$\Phi=uwva$. 
Since $uv$, $uw$, $vw$, $ua$, and $va$ are red edges,  
and $G$ is a simple graph, both $e$ and $f$ have $a$ as an endpoint, 
and then $w$ as the other endpoint, which is impossible.
This concludes the proof of Claim \ref{haromszog}.
\end{proof}

Return to the proof of Theorem \ref{thm:4n-8}. 
We apply a discharging argument. 
Give charge $2$ to every blue edge. 
Obviously, the sum of all charges is $2E_b$. We redistribute the charges in two steps. 
For each blue edge $e$, let $f$ be its odd pair. It is a red edge, 
let $\Phi_1$ and $\Phi_2$ be its incident faces ($\Phi_1$ and $\Phi_2$ might coincide).
Move charge $1$ to both $\Phi_1$ and $\Phi_2$ from $e$. 
For any face $\Phi$, let $ch(\Phi)$ denote the charge of $\Phi$.

Now for all triangular faces $\Delta$ 
with $ch(\Delta)>1$ do the following. 
By Claim \ref{haromszog}, there are two or three blue edges that belong to $\Delta$ and they all belong also to an adjacent face $\Phi$, which is not a triangle. Move $ch(\Delta)-1$ charge from $\Delta$ to $\Phi$. After this redistribution, 
for any face  $\Phi$, let $ch'(\Phi_i)$ denote its charge.
Clearly, the sum of all charges is still $2E_b$.
For any face $\Phi$, let $|\Phi|$ be the number of edges on its boundary, counted with multiplicity.

\begin{claim}\label{ch'}
For any face $\Phi$, $$ch'(\Phi)\le 3|\Phi|-8.$$
\end{claim}

\begin{proof}
Since $G$ is a simple graph, $|\Phi|\ge 3$. If $|\Phi|=3$, then trivially $ch'(\Phi)\le 1=3\cdot 3-8$.
Suppose that $|\Phi|\ge 4$. 
In the first redistribution, each unit of charge in $ch(\Phi)$ corresponds to a blue edge $e$ that belongs to 
$\Phi$, and its odd pair, $f$ is a boundary edge of $\Phi$. Therefore, $ch(\Phi)\le |\Phi|$.
If  $|\Phi|=4$, then 
by Claim \ref{haromszog}, $\Phi$ did not receive any charge in the second redistribution,
so $ch'(\Phi)=ch(\Phi)$. Consequently, $ch'(\Phi)=ch(\Phi)\le 4=3\cdot 4-8$.

Let $|\Phi|=5$, we have to prove that $ch'(\Phi)\le 7$. 
Then $ch(\Phi)\le 5$, so if it received at most $2$ charges in the second redistribution, we are done. 
Suppose that it received at least $3$ charges in the second round. Then it is adjacent to at least two triangles. By Claim \ref{haromszog}, if 
it receives one (two) charges from a triangle, then they have at least two (three) common sides. Therefore, 
the only possibility is that $\Phi$ is adjacent to two triangles, and one has three common sides with $\Phi$ and one has two. However, 
this is also impossible since $G$ is a simple graph.

Let $|\Phi|\ge 6$ now. Then $ch(\Phi)\le |\Phi|$, and it receives at most $2|\Phi|/3$ charges in the second round, so
$ch'(\Phi)\le 5|\Phi|/3\le 3|\Phi|-8$. 
\end{proof}

Using Claim \ref{ch'} we can obtain the desired bound on the number of edges.
$$E(G)\le E(H)+E_b\le E(H)+\frac{1}{2}\sum_{\mbox{\scriptsize {$\Phi$\ is a face}}} \left( 3|\Phi|-8\right) =
4E(H)-4F(H)=4n-8.$$

This 
proves Theorem 
\ref{thm:4n-8} if we assume that $G$ has a $1$-odd-planar drawing with a 
connected spanning even-crossing subgraph.

\smallskip

{\bf Assume now that $G$ does not have a $1$-odd-planar drawing with a 
connected spanning even-crossing subgraph.}

For any subgraph of $H\subset G$, an {\em external bridge} of $H$ is an edge of $G\setminus H$ between two connected components of $H$.
Take a 1-odd-plane drawing of $G$. Let $H$ be an even-crossing subgraph of $G$ with the minimum number of 
connected components. Among all such subgraphs, choose $H$ to have the maximum number of edges.

Observe that 
for every pair of edges $e,f$ forming an odd pair, exactly one of them is in $H$, 
otherwise $H+e$ (or $H+f$) would be an even-crossing subgraph with more edges.
For the same reason,  every edge of $G$ that does not have an odd pair is in $H$.

Define procedure {\sc Cut($G, H$)} recursively, as follows.
Let $G_0=G$, $H_0=H$. Suppose that we already have $G_i$, $H_i$, for all $i\ge 0$. 
If there is no external bridge of $H_i$ in $G_i$, then {\em Stop}. Otherwise, 
select an external bridge $e_i\in G_i$ of $H_i$ and 
let $f_i$ be its odd pair. Delete both of them, that is, let
$G_{i+1}=G_i-e_i-f_i$ and let 
$H_{i+1}=H_i-f_i$. 
Increase $i$ by $1$ and repeat.
Observe that it is not a deterministic procedure. 

\begin{claim}\label{oddpair}
In each step $i$, $f_i$, the odd pair of the external bridge $e_i$, is a cut-edge of $H_i$.  
\end{claim}

\begin{proof}
First we show that $f_i\in H_i$. 
Suppose that for some $i$, $e_i$ is an external bridge of $H_i$ in $G_i$, the edge $f_i$ is its odd pair, and both $f_i$ and $e_i$ are in $G_i \setminus H_i$.
Observe that throughout the procedure, edges are only deleted from the original graphs $G$ and $H$, while their colors remain unchanged. Consequently, $e_i$ and $f_i$ also form an odd pair in $G$, and neither belongs to the original subgraph $H$, a contradiction to the choice of $H$.

Now we show that $f_i$ is a cut-edge of $H_i$. Suppose that there is a counterexample, that is, for some $i$,
$f_i$ is not a cut-edge of $H_i$. So $f_i$ belongs to a cycle in $H_i$.
Take all possible procedures {\sc Cut($G, H$)} for every choice of even-crossing $H$ with the minimum number of components and among those, the maximum number of edges, which lead to a counterexample.
Take the choice and procedure which leads to a counterexample for the smallest value of $i$. 
With a slight abuse of notation, assume that this counterexample is $H=H_0, H_1, \ldots, H_i$ with deleted edges $e_0, e_1, \ldots e_i$ and $f_0, f_1, \ldots f_i$.

Suppose first that $i=0$, that is, $e_0$ is an external bridge of 
$H=H_0$ and its odd pair, $f_0\in H$, is not a cut-edge of $H$. 
Then let $H'=H+e_0-f_0$. 
The subgraph $H'$ is 
even-crossing, it has the same number of edges as $H$, and has fewer components,  contradicting the choice of $H$. So we have $i>0$. 

Let $j$ be the first index at which the endpoints of $e_i$ belong to different components of $H_j$. Clearly, $0\leq j\leq i$. 

Suppose that $j=0$. Let $H'=H+e_i-f_i$. Since $f_i$ belongs to a cycle in $H_i$, it also belongs to a cycle in $H$, therefore, $H-f_i$ has the same components as $H$. Since $j=0$, it implies that
$H'=H+e_i-f_i$ has fewer components than $H$. Since it is an even-crossing subgraph, it contradicts the choice of $H$. So, we have $j>0$ and the subgraph 
$H'$ has the same components as $H$, the same number of edges, and $H'$ is also an even-crossing subgraph.
 
Consider  procedure {\sc Cut($G, H'$)} with the same choices of 
edges  $e_0, e_1, \ldots e_{j-1}$ and $f_0, f_1, \ldots f_{j-1}$. This is possible by the assumptions. Let $H'=H_0', H_1', \ldots, H_{j-1}'$ be the resulting subgraphs.

We compare the two procedures {\sc Cut($G, H$)} and {\sc Cut($G, H'$)}.

In {\sc Cut($G, H$)},  
at step $j-1$, the external bridge $e_{j-1}$ has odd pair $f_{j-1}$, which is a cut-edge of $H_{j-1}$ by the minimality of $i$. By the definition of $j$, the endpoints of $e_i$ are connected in $H_{j-1}$ but disconnected in $H_j=H_{j-1}-f_{j-1}$. 

Therefore, they belong to the two components created by deleting $f_{j-1}$. Adding $e_i$ reconnects these components. Moreover, since $f_i$ belongs to a cycle of $H_j$, deleting it does not change the connected components. Consequently, $f_{j-1}$ is not a cut-edge of $H'_{j-1}=H_{j-1}+e_i-f_i$. Thus, we obtain a counterexample in {\sc Cut($G, H'$)},
at step $j-1<i$, contradicting the minimality of $i$.
This proves the claim.

\end{proof}

Now we can finish the proof of Theorem \ref{thm:4n-8}.
Apply procedure {\sc Cut($G$, $H$)}. Suppose that it stops at step $i$.
Suppose that $H_i$ has $a$ isolated vertices, $b$ components with two vertices, and $c$ components
with at least $3$ vertices. Since $H$ is not connected and in each step we increase the number of components by one by Claim \ref{oddpair}, 
we have $a+b+c\ge i+2$. 
Let $U$ be the vertex set of a connected component of $H_i$. 
If $|U|=1$ then clearly $|E(G(U))|=0$, and if $|U|=2$ then $|E(G(U))|=1$. 
If $|U|\ge 3$, then $G(U)$ is a 1-odd-planar graph with a connected even-crossing subgraph. Therefore, by the first part of the proof, we have $|E(G(U))|\le 4|U|-8$. 
Let $V_1, \ldots, V_c$ be the vertex sets of the components of at least $3$ vertices. 
We have $\sum_1^c|V_i|=n-a-2b$.
 
In each step of {\sc Cut($G$, $H$)} we deleted two edges, so $|E(G_i)|=|E(G)|-2i$.
And since the procedure stopped at step $i$, there are no edges in $G_i$ between the components of $H_i$. 
Summing it up, we have
$$|E(G)|\le 2i+|E(G_i)|\le 2(a+b+c-2)+b+\sum_{j=1}^c\left(4|V_j|-8\right)$$
$$=
                           2a+2b+2c-4+b+4(n-a-2b)-8c=
4n-2a-5b-6c-4.$$

It remains to show that $4n-2a-5b-6c-4\le 4n-8$. 
It holds trivially if $a\ge 2$ or
$b\ge 1$ or $c\ge 1$ and at least one of these conditions should hold, since $n\ge 3$.
This finishes the proof of Theorem \ref{thm:4n-8}.
\qed

\section{Proof of Theorem \ref{oddcrlemma}}\label{sec:crlemma}

The following proof is analogous to that of Theorem 1 in \cite{KT23}, but Theorem \ref{thm:4n-8} allows us to obtain a stronger result.

\begin{proof}
    We have $m_0^{\mathrm{odd}}(n) = 3n-6$, and by Theorem \ref{thm:4n-8}, $m_1^{\mathrm{odd}}(n) = 4n-8$.
    We first prove by induction on the number of edges that every graph $G$ with $n$ vertices and $m$ edges satisfies $\ocr(G)\ge m-3n$. If $m\le 3n$, the claim is trivial. Suppose that $m>3n$ and that the claim holds for graphs with $n$ vertices and $m-1$ edges. Consider an arbitrary drawing $D$ of $G$. Since $m_0^{\mathrm{odd}}(n)=3n-6$, the drawing $D$ contains an odd pair of edges $e$ and $f$. Remove $e$ from $G$ and consider the resulting drawing. It has one fewer edges and at least one fewer odd pair. By the induction hypothesis, the number of odd pairs in $D$ is therefore at least $1+(m-1)-3n=m-3n$.

    We next prove, again by induction on the number of edges, that every graph $G$ with $n$ vertices and $m$ edges satisfies $\ocr(G)\ge 2m-7n$. If $m\le 4n$, then the first bound gives $\ocr(G)\ge m-3n\ge 2m-7n.$
    Suppose that $m>4n$ and that the claim holds for graphs with $n$ vertices and $m-1$ edges. Consider an arbitrary drawing $D$ of $G$. Since $m_1^{\mathrm{odd}}(n)=4n-8$, there is an edge $e$ that belongs to two odd pairs, namely ${e,f}$ and ${e,g}$. Remove $e$ from $G$ and consider the resulting drawing. It has one fewer edge and at least two fewer odd pairs. By the induction hypothesis, the number of odd pairs in $D$ is therefore at least $2+2(m-1)-7n=2m-7n$.

    Let $G$ be a graph with $n$ vertices and $m\ge 5.25n$ edges, drawn in the plane so that every pair of adjacent edges crosses an even number of times and exactly $\ocr_*(G)$ pairs of edges cross an odd number of times.
    Construct a random subgraph $G'$ by selecting each vertex independently with probability $p=5.25n/m$. Let $n'$, $m'$, and $x(G')$ denote, respectively, the number of vertices and edges of $G'$, and the number of odd-crossing pairs of edges in $G'$, in the inherited drawing.
    We have 
    $$E(n')=pn, \ \ E(m')=p^2m, \ \ E(\ocr_{*}(G'))\le E(x(G'))=p^4\ocr_{*}(G).$$
    For $G'$ we have  $\ocr_{*}(G')\ge \ocr(G')\ge 2m'-7n'$, taking expected values, 
    $$p^4\ocr_{*}(G)\ge 2p^2m-7pn.$$
    For  $p=5.25n/m$, we have 
    $$\ocr_{*}(G)\ge\frac{1}{41.34375}\frac{m^3}{n^2}.$$

\end{proof}

\bigskip

\noindent {\bf Remarks.}

1. We have seen that 
$m^{\mathrm{odd}}_0(n)=m_0(n)$ and 
$m^{\mathrm{odd}}_1(n)=m_1(n)$. In general, we know that $m^{\mathrm{odd}}_k(n)\ge m_k(n)$. 

\begin{problem}
    Does $m^{\mathrm{odd}}_k(n)= m_k(n)$ hold for every $k\ge 0$?
\end{problem}

The crossing number and the odd-crossing number are not always equal;  there is an infinite family of graphs with $\ocr(G) < \cro(G)$ \cite{PSS08}, \cite{T08}.
Therefore, this problem might be hard, and it is not easy even to have a conjecture on it. 

\smallskip

2. As we mentioned in the Introduction, 
the proof of the Crossing Lemma by 
Chazelle, Sharir and Welzl \cite{AZ99} works for all versions in Table 1,
but the improvements do not.
The reason is that for the probabilistic arguments, 
similar to the proof of Theorem \ref{oddcrlemma}, 
we need many crossings between {\em independent edges}. It is not clear, if it is required only for this 
argument, or is it really an essential requirement. The original argument can be saved for all versions, but the improvements it works only if we {\em do not allow} adjacent edges to cross
(or to cross an odd number of times). That's why we have
improvements only for $\pcr_{+}$, $\ocr_{+}$, and $\ocr_{*}$.

\begin{problem}
    Improve the constant in the Crossing Lemma for the other versions of the crossing number.
\end{problem}

\smallskip 

3. All results presented in this paper were obtained by the authors without the use of artificial intelligence tools.



\begin{thebibliography}{99}


\bibitem[A19]{A19}
E. Ackerman,
On topological graphs with at most four crossings per edge,
{\it Computational Geometry} {\bf 85} (2019), 101574.

\bibitem[AS14]{AS14}
E. Ackerman, M. Schaefer,
A crossing lemma for the pair-crossing number,
in: {\it Graph Drawing},
C. Duncan, A. Symvonis (eds.),
Lecture Notes in Computer Science {\bf 8871},
Springer, Berlin, Heidelberg, (2014), 222--233.

\bibitem[ACNS82]{ACNS82}
M. Ajtai, V. Chv\'atal, M. Newborn, E. Szemer\'edi,
Crossing-free subgraphs,
in: {\it Theory and Practice of Combinatorics},
North-Holland Mathematics Studies {\bf 60},
North-Holland, Amsterdam--New York, (1982), 9--12.

\bibitem[AZ99]{AZ99}
M. Aigner, G. Ziegler,
{\it Proofs from the Book},
Springer, Berlin, (1999).

\bibitem[BK26]{BK26}
A. B\"ungener, M. Kaufmann,
Improving the crossing lemma by characterizing dense
2-planar and 3-planar graphs,
{\it Journal of Graph Algorithms and Applications}
{\bf 29}(3) (2026), 143--174.

\bibitem[BFKP26]{BFKP26}
A. B\"ungener, J. Franz, M. Kaufmann, M. Pfister,
A first view on the density of 5-planar graphs,
arXiv preprint arXiv:2505.24364 (2026).

\bibitem[CSSW20]{CSSW20}
\'E. Czabarka, I. Singgih, L. Sz\'ekely, Z. Wang,
Some remarks on the midrange crossing constant,
{\it Studia Scientiarum Mathematicarum Hungarica}
{\bf 57}(2) (2020), 187--192.

\bibitem[H34]{H34}
C. Chojnacki (H. Hanani),
\"Uber wesentlich unpl\"attbare Kurven
im drei-dimensionalen Raume,
{\it Fundamenta Mathematicae} {\bf 23} (1934), 135--142.


\bibitem[F26]{F26}
R. Fulek,
Why do adjacent crossings matter?,
in: {\it New Probes into Discrete and Convex Geometry},
J. Pach, G. T\'oth (eds.),
Bolyai Society Mathematical Studies {\bf 33},
to appear.

\bibitem[KT23]{KT23}
J. Karl, G. T\'oth,
Crossing lemma for the odd-crossing number,
{\it Computational Geometry} {\bf 108} (2023), 101901.

\bibitem[L83]{L83}
F. T. Leighton,
{\it Complexity Issues in VLSI},
MIT Press, Cambridge, (1983).

\bibitem[PT97]{PT97}
J. Pach, G. T\'oth,
Graphs drawn with few crossings per edge,
{\it Combinatorica} {\bf 17}(3) (1997), 427--439.

\bibitem[PT00a]{PT00a}
J. Pach, G. T\'oth,
Which crossing number is it anyway?,
{\it Journal of Combinatorial Theory, Series B}
{\bf 80}(2) (2000), 225--246.

\bibitem[PT00b]{PT00b}
J. Pach, G. T\'oth,
Thirteen problems on crossing numbers,
{\it Geombinatorics} {\bf 9}(4) (2000), 194--207.

\bibitem[PRTT06]{PRTT06}
J. Pach, R. Radoi\v{c}i\'c, G. Tardos, G. T\'oth,
Improving the crossing lemma by finding more crossings
in sparse graphs,
{\it Discrete and Computational Geometry}
{\bf 36}(4) (2006), 527--552.

\bibitem[PSS08]{PSS08}
M. Pelsmajer, M. Schaefer, D. \v{S}tefankovi\v{c},
Odd crossing number and crossing number are not the same,
{\it Discrete and Computational Geometry}
{\bf 39} (2008), 442--454.

\bibitem[PSS07]{PSS07}
M. Pelsmajer, M. Schaefer, D. \v{S}tefankovi\v{c},
Removing even crossings,
{\it Journal of Combinatorial Theory, Series B}
{\bf 97}(4) (2007), 489--500.

\bibitem[S12]{S12}
M. Schaefer,
The graph crossing number and its variants: A survey,
{\it Electronic Journal of Combinatorics}
(2012), DS21.

\bibitem[S18]{S18}
M. Schaefer,
{\it Crossing Numbers of Graphs},
CRC Press, (2018).

\bibitem[T08]{T08}
G. T\'oth,
Note on the pair-crossing number and the odd-crossing number, 
{\it Discrete and Computational Geometry} {\bf 39.4} (2008), 791--799.

\bibitem[T26]{T26}
G. T\'oth,
Generalizations of the crossing lemma,
in: {\it Sum(m)it280: Surveys in Extremal
Combinatorics and Combinatorial Geometry},
G. O. H. Katona, B. Patk\'os, C. Tompkins (eds.),
Bolyai Society Mathematical Studies {\bf 32},
Springer Nature Switzerland, Cham, (2026), 385--411.

\end{thebibliography}
\end{document}